\documentclass[11pt]{article}
\usepackage{amssymb,amsthm,amscd,array}
\usepackage{graphics}
\usepackage[final]{epsfig}
\usepackage{amsmath, amssymb}
\usepackage{amsthm}
\usepackage{pb-diagram}
\usepackage{amsmath,amsthm,amsfonts,amssymb,latexsym}
\usepackage{a4wide}

\def\dis{\displaystyle}
\numberwithin{equation}{section}

\newtheorem{theorem}{Theorem}[section]
\newtheorem{definition}[theorem]{Definition}

\newtheorem{remark}[theorem]{Remark}
\newtheorem{proposition}[theorem]{Proposition}

\newcommand{\R}{\mathbb{R}}

\newcommand{\N}{\mathbb{N}}

\newcommand{\C}{\mathbb{C}}

\newcommand{\fL}{{\mathfrak{L}}}

\newcommand{\cK}{{\mathcal{K}}}

\newcommand{\Hom}{\mathrm{Hom}}

\def\d{\delta}

\def\a{\alpha}

\def\s{\sigma}

\def\l{\lambda}
\def\m{\mu}

\title{Conformally equivariant quantization and symbol  maps
associated with $n$-ary differential operators on weighted densities
  } \vspace{8pt}
\author{T. Bichr
\thanks{
D\'epartement de Math\'ematiques, Facult\'e des sciences de Sfax,
3000 Sfax BP 1171, taher-bechr@hotmail.fr} \and J. Boujelben
\thanks{
D\'epartement de Math\'ematiques, Facult\'e des sciences de Sfax,
3000 Sfax BP 1171, Tunisie, jamel\_boujelben@hotmail.fr} \and
 K.Tounsi
\thanks{
D\'epartement de Math\'ematiques, Facult\'e des sciences de Sfax,
3000 Sfax BP 1171, Tunisie, Khaled\_286@yahoo.fr }}
\date{}

\begin{document}\maketitle
\textbf{Abstract} We are interested in the study of the space of
$n$-ary differential operators denoted by
$\mathfrak{D}_{\underline{\l},\mu}$ where
$\underline{\l}=(\l_{1},...,\l_{n})$ acting on weighted densities
from $\frak F_{\l_1}\otimes\frak F_{\l_2}\otimes...\otimes\frak
F_{\l_n}$ to $\frak F_{\mu}$ as a module over the orthosymplectic
superalgebra $\mathfrak{osp}(1|2)$. As a consequence, we prove the
existence and the uniqueness of a canonical conformally equivariant
symbol map from $\mathfrak{D}_{\underline{\lambda},\mu}^k$ to the
corresponding space of symbols as well for the  explicit expression
of the associated quantization map.

\vspace{15pt}



\section{Introduction}
The quantization is a concept that comes from physics. The quantization of a classical system
whose phase space is a symplectic manifold, consists in the construction of a Hilbert
space H and a correspondence between classical and quantum observables.
 Let $M$ be a smooth manifold, $T^*M$  the cotangent bundle on
 $M$ and $\mathcal{S}(M)$  the space of smooth functions on $T^*M$ polynomial on the fibers, which the  is usually called the space of symbols of differential
operators. The standard quantization procedure consists of
constructing  a map $\mathcal{Q}$ between the space
$\hbox{Pol}(T^*M)$ of polynomials on $T^*M$ and the space
$\mathcal{D}(M)$ of linear differential operators on $M$ called a
\textit{quatization map}. The inverse $\sigma= \mathcal{Q}^{-1}$ is
thus called a \textit{symbol map}. Generally, there is no
quantization  and symbol maps equivariant with respect to the action
 the Lie algebra $\mathrm{Vect}(M)$ of vector fields on $M$ (or the group $\mathrm{Diff}(M)$
 of diffeomorphisms of $M$) on the two spaces $\mathcal{D}(M)$ and
 $\hbox{Pol}(T^*M)$. Thus, we restrict ourselves to equivariant symbols and
quantization maps with respect to the action of a given  subalgebra
of $\mathrm{Vect}(M)$.

More precisely, Let
 for every $\l\in \C$,  ${\cal F}_\lambda(M)$ and ${\cal D}_{\lambda,\mu}(M)$ stand for the space
of tensor densities of degree $\l$ on $M$ and   the space of linear
differential operators from ${\cal F}_\lambda(M)$ to ${\cal
F}_\m(M)$ ($\lambda, \mu\in \C$) respectively.
 These spaces are naturally modules over the Lie algebra $\mathrm{Vect}(M)$.The space of
symbols corresponding to ${\cal D}_{\lambda,\mu}(M)$ is there for $\mathcal{S}_{\delta}(M)=\mathcal{S}(M)\otimes{\cal F}_\delta(M)$ where $\delta=\mu-\lambda$,there
 is a filtration
 \begin{equation*}
{\cal D}^{0}_{\lambda,\mu}\subset{\cal D}^{1}_{\lambda,\mu}\subset\cdots
{\cal D}^{k}_{\lambda,\mu}\subset \cdots.
\end{equation*}
and the associated module $\mathcal{S}_{\delta}(M)=gr({\cal D}_{\lambda,\mu})$ is graded by the degree of polynomials:
\begin{equation*}
\mathcal{S}^{0}_{\delta}\subset\mathcal{S}^{1}_{\delta}\subset\cdots
\mathcal{S}^{k}_{\delta}\subset \cdots.
\end{equation*}
The problem of equivariant quantization is the quest for a
quantization map:
$$Q_{\lambda,\mu}:\mathcal{S}_{\delta}(M)\longrightarrow{\cal D}_{\lambda,\mu}(M)$$
that commutes with the action of a given lie subalgebra of
$\mathrm{Vect}(M)$. In other word, it amounts to an identification
of these two spaces which is canonical with the respect to the
geometric on $M$. The inverse of the quantization map.
$$ \sigma_{\lambda,\mu}:=(Q_{\lambda,\mu})^{-1}$$
is called symbol map.

The concept of equivariant quantization over $\R^n$ was introduced
by P. Lecomte and V. Ovsienko in \cite{LO}. In this seminal work,
they considered spaces of differential operators acting between
densities and the Lie algebra of projective vector fields over
$\R^n$, $\mathrm{sl}(n + 1)$. In this situation, they showed the
existence and uniqueness of an equivariant quantization. This
results were generalized in many references (see for instance
\cite{DLO}, \cite{L1}). In \cite{L2}, P. Lecomte globalized the
problem of equivariant quantization by defining the problem of
natural invariant quantization on arbitrary manifolds. Finally in
\cite{BBT}, \cite{Bo}, \cite{Bos}, \cite{CS}, \cite{F}, \cite{H1},
\cite{H2}, \cite{MR1}, \cite{MR2}, \cite{MR3}, \cite{MR4}, the
authors proved the existence of such quantizations by using
different methods in more and more general contexts.

Recently, several papers dealt with the problem of equivariant
quantizations in the context of supergeometry: the papers \cite{LMR}
and \cite{MR5}  exposed and solved respectively the problems of the
$\frak {pgl}(p + 1|q)$-equivariant quantization over the superspace
$\R^{p|q}$ and of the $\frak{osp}(p + 1; q + 1|2r)$-equivariant
quantization over $\R^{p+q|2r}$, whereas in \cite{LR}, the authors
define the problem of the natural and projectively invariant
quantization on arbitrary supermanifolds and show the existence of
such a map. In \cite{GMO}, \cite{Me}, \cite{MNR} the problem of
equivariant quantizations over the supercircles $S^{1|1}$ and
$S^{1|2}$ endowed with canonical contact structures was considered,
these quantizations are equivariant with respect to Lie
superalgebras $\mathfrak{osp}(1|2)$ and $\mathfrak{osp}(2|2)$ of
contact projective vector fields respectively.

In \cite{BBTT}, for the $S^{1|1}$-case, we were interested in the
study of the space $\frak {D}_{\lambda_1,\lambda_2,\mu}$ of bilinear
differential operators from $\frak F_{\l_1}\otimes\frak F_{\l_2}$ to
$\frak F_{\mu}$. For almost all values $(\l_1,\l_2,\mu)$, we prove
the existence and the uniqueness (up to normalization) of a
projectively ,i.e., $\mathfrak{osp}(1|2)$-equivariant symbol map
between $\frak {D}_{\lambda_1,\lambda_2,\mu}$ and the corresponding
space of symbols $\frak {S}_{\lambda_1,\lambda_2,\mu}$   and
calculate the explicit expressions  of the symbol and the associated
quantization maps.

Our motivation  in this work  is the generalization  of the results
proved \cite{BBTT}. Namely we consider the
 superspace $\frak{D}_{\underline{\lambda},\mu}$, $\underline{\lambda}=(\lambda_1,\cdots, \lambda_n)\in \C^n$, of $n$-ary differential operators  $A:\frak F_{\l_1}\otimes\frak F_{\l_2}\otimes...\otimes\frak F_{\l_n}\rightarrow\mathfrak{F}_\mu$,
 where $\frak F_{\l}, \l\in \C$, is the space of tensor densities on the supercircle $S^{1|1}$ of degree $\lambda$. The analogue, in the super setting, of
the projective algebra $\hbox{sl}(2)$ is  the orthosymplectic Lie
superalgebra $\mathfrak{osp}(1|2)$, which is the smallest simple Lie
superalgebra, can be realized as a subalgebra of
$\mathrm{Vect}_\C(S^{1|1})$. Naturally, the Lie superalgebra
 $\mathrm{Vect}_\C(S^{1|1})$, and therefor $\mathfrak{osp}(1|2)$, act on $\frak{D}_{\underline{\lambda},\mu}$, the $\mathfrak{osp}(1|2)$-module
$\frak{D}_{\underline{\lambda},\mu}$ is filtered as:
\begin{equation*}
\mathfrak{D}^{0}_{\underline{\lambda},\mu}\subset\mathfrak{D}^{\frac{1}{2}}_{\underline{\lambda},\mu}\subset
\mathfrak{D}^{1}_{\underline{\lambda},\mu}\subset
\mathfrak{D}^{\frac{3}{2}}_{\underline{\lambda},\mu}
\subset\cdots\subset\mathfrak{D}^{k-\frac{1}{2}}_{\underline{\lambda},\mu}\subset
\mathfrak{D}^{k}_{\underline{\lambda},\mu}\subset \cdots.
\end{equation*}
The graded module
$\hbox{gr}(\mathfrak{D}_{\underline{\lambda},\mu})$, also called the
space of symbols and denoted by
$\mathcal{S}_{\underline{\lambda},\mu}$, depends only on the shift,
$\d=\mu- |\underline{\l}|$, $|\underline{\l}|= \lambda_1+\cdots
+\lambda_n$, of the weights. Moreover,  as a
$\mathrm{Vect}_\C(S^{1|1})$-module,
$\mathcal{S}_{\underline{\lambda},\mu}$ is decomposed as
$\dis\bigoplus_{k\in \frac{1}{2}\N}
\mathcal{S}_{\underline{\lambda},\mu}^k$ where $${\cal
S}^k_{\underline{\lambda},\mu}=\bigoplus_{\ell=0}^{2k}\mathfrak{D}^{\ell}_{\underline{\lambda},\mu}/\mathfrak{D}^{\ell-\frac{1}{2}}_{\underline{\lambda},\mu}=\bigoplus_{\ell=0}^{2k}
{\mathfrak{F}}^{(\ell)}_{\delta-\frac{\ell}{2}},$$
  ${\mathfrak{F}}^{(\ell)}_{\delta-\frac{\ell}{2}}$  stands for the sum $ \bigoplus\mathfrak{F}_{\delta-\frac{\ell}{2}}$ where
$\mathfrak{F}_{\delta-\frac{\ell}{2}}$ is counted  $\left(\ell +n-1 \atop n-1  \right)$
times.\\
Moreover, we prove that,  if $\delta=\mu-|\underline{\l}|\neq
\dis\frac{1}{2}, 1,\frac{3}{2},2,\frac{5}{2},\cdots,k$, then
$\mathfrak{D}^k_{\underline{\l},\mu}$ is isomorphic to
$\mathcal{S}_{\underline{\l},\mu}^k$  as an
$\hbox{osp}(1|2)$-module. This isomorphism, called a \textit{a
conformally equivariant symbol map}, is unique (once we fix a
principal symbol). Explicit expressions of the normalized symbol and
its inverse, the \textit{ conformally equivariant quatization map},
are given.

\section{The main definitions }
In this section,we recall the main definition and facts related to
the geometry of the supercircle $S^{1|1}$.(See for instance
\cite{BBB}, \cite{GO}, \cite{GMO})

\subsection{Geometry of the supercircle $S^{1|1}$}
 The supercircle $S^{1|1}$ is the
simplest supermanifold of dimension $1|1$ generalizing $S^1$. In order to fixe notation,
let us give here  the basic definitions of geometric objects on  $S^{1|1}$ . We define the supercircle $S^{1|1}$  by describing its graded commutative algebra of functions which
we denote by $C^\infty_\C (S^{1|1})$ and which is constituted by the elements
\begin{equation}
F= f_0(x)+\theta f_1(x), \end{equation} where
$x$ is an arbitrary parameter on $S^1$ (the even variable), $\theta$
is the odd variable ($\theta^2=0$) and $f_0$, $f_1$ are $C^\infty$
complex valued functions. We denote by $F'$ the derivative of $F$
with respect to $x$, i.e, $F' (x,\theta)=
f'_0(x)+\theta f'_1(x) $.

\subsection{Vector fields and differential forms}

Let $\mathrm{Vect}({S}^{1|1})$ be the
superspace of vector fields on ${S}^{1|1}$:
\begin{equation}
\mathrm{Vect}_\C(S^{1|1})=\left\{F_0\partial_x+ F_1\partial_\theta
\mid ~F_i\in C^\infty_\C(S^{1|1})\right\}, \end{equation} where
$\partial_\theta$ (resp $\partial_x$) means the partial derivative
$\dis\frac{\partial}{\partial\theta}$ (resp
$\dis\frac{\partial}{\partial x}$).\\
Let $\Omega ^{1}({S}^{1|1})$be the rank $1|1$ right
$C^\infty_\C(S^{1|1})$-module
 with basis $dx$ and $d\theta$, we interpret it as the
right dual over
$C^\infty_\C(S^{1|1})$ to the left $C^\infty_\C(S^{1|1})$-module
 $\mathrm{Vect}_\C(S^{1|1})$, by setting \\ $\langle\partial_{y_{i}}, dy_{j}\rangle=\delta_{ij}$  for
$y = (x, \theta)$. The space
  $\Omega ^{1}({S}^{1|1})$  is a left module over $\mathrm{Vect}_\C(S^{1|1})$, the action being given by the
Lie derivative:
 $$ \langle X,L_{Y}(\alpha)\rangle:=\langle [X,Y],\alpha\rangle  $$

\subsection{Lie superalgebra of contact vector fields}
The standard contact structure on $S^{1|1}$ is defined as a codimension 1 non-integrable distribution
$\langle\overline{D}\rangle$ on $S^{1|1}$, i.e., a subbundle in $TS^{1|1}$ generated by the odd vector field
\begin{equation}  \overline{D}= {\partial_\theta} - \theta{\partial_ x},\end{equation}
This contact structure can be equivalently defined as the kernel of
the differential 1-form
\begin{equation}\a=dx+\theta d\theta.\end{equation}
These vector fields satisfy the condition \begin{equation}
\overline{D}^{2j}= (-1)^jD^{2j}=(-1)^j
\partial_x^j ,  \forall j\in \N.\end{equation}
where $D= {\partial_\theta} + \theta{\partial_ x}$.\\
One can easily check the \textit{super Leibniz formula}:
\begin{equation}\label{super Leibniz formula}\overline{D}^j \circ F = \dis\sum _{i=0}^j \left( j \atop i
\right)_s
(-1)^{|F|(j-i)}\overline{D}^i(F)\overline{D}^{j-i},\end{equation}
where  the notions $\left( j \atop i  \right)_s$ and $|~|$ stand
respectively  for the \textit{super combination} defined by
\begin{equation}\label{super
combination}\left( j \atop i \right)_s = \left\{
                                  \begin{array}{ll}
                                     \left( [\frac{j}{2}] \atop [\frac{i}{2}]  \right)& \hbox{if i is even or j is odd} \\
                                    0  & \hbox{otherwise.}
                                  \end{array}
                                \right.\end{equation}
and  for the parity function ($[x]$ denotes the integer part
of a real number $x$). \\
 A vector field $X$ is said to be contact
if it preserves the contact distribution, i.e.,
\begin{equation}
[X,\overline{D}]=
F_X\overline{D},\end{equation}
where $F_X \in C^\infty_\C(S^{1|1})$ is a function depending on X.\\
We denote by $\mathcal{K}(1)$ the \textit{Lie superalgebra of
contact vector fields} on  $S^{1|1}$.
It is well-known that every contact vector field can be expressed, for some function
$f\in C^\infty_\C(S^{1|1})$, by (see \cite{GMO}):
 \begin{equation}\label{act} X_f=
-f\overline{D}^2 + \frac{1}{2}D(f)\overline{D} .\end{equation}  The
vector field (\ref{act}) is said to be the contact vector field with
contact Hamiltonian f. One checks that
$$L_{X_{f}}\a=f'\a~~,~~~~[X_{f},\overline{D}]=-\frac{1}{2}f'\overline{D}$$
The contact bracket is defined by $[X_{f} ,X_{g}] = X_{\{f,g\}}$. The space $C^\infty_\C(S^{1|1})$ is thus equipped
with a Lie superalgebra structure isomorphic to $\mathcal{K}(1)$. The explicit formula can be easily calculated:
\begin{equation}\{f,g\}=fg'-f'g+\frac{1}{2}(-1)^{|f|(|g|+1)}D(f)D(g).\end{equation} The action of $\mathcal{K}(1)$ on
$C^{\infty}_\C(S^{1|1})$ is defined by: \begin{equation} \frak
L_{X_f}(g)=f g'+ \frac{1}{2}D(f)\overline{D}(g) =f
g'+\frac{1}{2}(-1)^{|f|+1}\overline{D}(f)\overline{D}(g).
\end{equation}
\subsection{The orthosymplectic Lie  superalgebra
$\mathfrak{osp}(1|2)$} If we identify $S^1$ with $\mathbb{RP}^1$
with homogeneous coordinates $(x_1:x_2)$ and choose the affine
coordinate $x= x_1/x_2$, the vector fields
$$\dis\frac{d}{dx},\,x\frac{d}{dx},\,x^2\frac{d}{dx}$$ are globally
defined and correspond to The standard projective structure on
$\mathbb{RP}^1$. In this adapted coordinate the action of the
subalgebra $\frak {sl}(2)$ of the Lie algebra ${{\rm Vect}}(S^1)$:
\begin{equation*}\frak {sl}(2)=\text{Span}\left(\frac{d}{dx},\,x\frac{d}{dx},\,x^2\frac{d}{dx}\right)\end{equation*}
is well defined.\\
 Similarly, we consider the orthosymplectic
Lie superalgebra $\mathfrak{osp}(1|2)$ as a subalgebra of
$\mathcal{K}(1)$:
\begin{equation}
\mathfrak{osp}(1|2)=\text{Span}(X_1=\partial_x,\,X_{x}= x\partial_x +
\frac{1}{2}\theta\partial_\theta
,\,X_{x^2}=x^2\partial_x +
x\theta\partial_\theta,\, X_{\theta}=
\frac{1}{2}D,\,X_{x\theta}= \frac{1}{2}xD). \end{equation} The space
of even elements :
\begin{equation}
(\mathfrak{osp}(1|2))_0=\text{Span}(X_1,\,X_{x},\,X_{x^2})
\end{equation} is isomorphic to $\mathfrak{sl}(2)$,
the space of odd elements is two dimensional:
\begin{equation}
(\mathfrak{osp}(1|2))_1=\text{Span}(
X_{\theta}=D,\,X_{x\theta}=xD).
\end{equation}
The new commutation relations are
\begin{equation*}\begin{array}{llll}
&[X_{1}, X_{x^2}]=2X_{x},~~&[X_\theta, X_\theta]=\frac{1}{2}X_{1},~~&[X_{x},X_1]=-X_{1},\\
&[X_x, X_{x^2}]= X_{x^2},~~&[X_{x\theta},X_{x\theta}]= \frac{1}{2}X_{x^2},~~&[X_{x^2},X_\theta]= -X_{x\theta},\\
&[X_x,X_\theta]=-\frac{1}{2}
X_\theta,~~&[X_1,X_{x\theta}]= X_\theta,~~&[X_1,X_\theta]=0,\\
&[X_x,X_{x\theta}]=\frac{1}{2}
X_{x\theta},~~&[X_{x^2},X_{x\theta}]=0,~~&[X_{x\theta},X_{\theta}]= \frac{1}{2} X_{x}.
\end{array}
\end{equation*}
 As in the $S^1$ case, there exist
adapted coordinates $(x,\theta)$ for which the
$\mathfrak{osp}(1|2)$-action is well defined (see \cite{GMO} for more
details).
\subsection{The space of weighted densities on $S^{1|1}$}
In the super setting, by replacing $dx$ by the 1-form $\alpha$, we
get analogous definition for weighted densities, {\sl i.e.}, we
define the space of $\lambda$-densities as
\begin{equation}
\mathfrak{F}_\lambda=\left\{F\alpha^\lambda~~|~~F \in
C^\infty_\C(S^{1|1})\right\}.
\end{equation}
As a vector space, $\mathfrak{F}_\lambda$ is isomorphic to
$C^\infty_\C(S^{1|1})$.\\
 For contact vector field $X_F$, define a
one-parameter family of first order differential operator on
$C^\infty_\C(S^{1|1})$

\begin{equation}\label{actiondensite}
\mathfrak{L}^{\lambda}_{X_F}= X_F + \lambda F', \lambda \in \C.
\end{equation} One easily checks that the map
$X_F \mapsto \mathfrak{L}^{\lambda}_{X_F}$ is a homomorphism of Lie
superalgebra , {\sl i.e.}, $
[\mathfrak{L}^{\lambda}_{X_F},\mathfrak{L}^{\lambda}_{X_G}] =
\mathfrak{L}^{\lambda}_{[{X_F},\,X_G]}$, for every $\lambda$. Thus
$\mathfrak{F}_\lambda$ becomes a $\mathcal{K}(1)$-module on
$C^\infty_\C(S^{1|1})$. Evidently, the Lie derivative of the density
$G\alpha^\lambda$ along the vector field $X_F$ in $\mathcal{K}(1)$
is given by:
\begin{equation}
\mathfrak{L}^{\lambda}_{X_F}(G\alpha^\lambda)=(X_F(G)+ \lambda
F'G)\alpha^\lambda. \end{equation} Explicitly, if we put
$F=f_0(x)+f_1(x)\theta$, $G=g_0(x)+g_1(x)\theta$,
\begin{equation}
\mathfrak{L}^{\lambda}_{X_F}(G)=L^{\lambda}_{f_0\partial_x}(g_0)+\frac{1}{2}~f_1g_1+
\left(L^{\lambda+\frac{1}{2}}_{f_0\partial_x}(g_1)+\lambda g_0f_1'+\frac{1}{2} g'_0 f_1\right)\theta.
\end{equation}

\subsection{Multilinear differential operators on weighted densities}
 We fix a natural number $n$. In order to avoid clutter, we have
found that it is convenient to use the notations of \cite{Bo}:
\begin{itemize}
  \item Denote by i either the n-tuple $(i_{1}, \cdots , i_{n})$ or the indices $i_{1}, \cdots , i_{n}$, as, for instance,\\
$a_{\underline{i}} = a_{i_{1}, \cdots , i_{n}}$. The difference should be discernable from the context.
  \item Denote by $|\underline{i}|$ the sum $\sum _{j=1}^{n}i_{j}.$
  \item Denote $\mathbf{1}_{i} := (0, \cdots , 0, 1, 0,\cdots , 0)$, where 1 is in the i-th position.
  \item Denote by $ \mathfrak{S}^{(i)}_{\l}= \oplus \mathfrak{F}_{\l}$ where $\mathfrak{F}_{\l}$ is
  counted $\left( i+n-1 \atop n-1  \right)$ times.
  \item $\dis\bigotimes_{t=1}^{n}\overline{D}^{i_{t}} :=\overline{D}^{i_{1}}\otimes\overline{D}^{i_{2}}\otimes\cdots\otimes\overline{D}^{i_{n}}. $
\item Throughout the text, we use  the classical convention  $\sum_{i=1}^0 c_i = 0$.
\end{itemize}

Obviously, $\forall \l_1,\l_2,\cdots,\l_n \in \R$,
$\mathfrak{F}_{\l_1}\otimes \mathfrak{F}_{\l_2}\otimes \cdots
\otimes\mathfrak{F}_{\l_n}$  also a $\mathcal{K}(1)$-module with the
action
\begin{equation} \label{Action on densities}\mathfrak{L}_{X_F}^{\underline{\l}}(\Phi_1\otimes\Phi_2\otimes...\otimes\Phi_n)=
\sum_{p=1}^{n}(-1)^{|F|(\sum_{i=1}^{p-1}|\Phi_{i}|)}
\Phi_1\otimes\Phi_2\otimes...\otimes\mathfrak{L}_{X_F}^{\l_p}(\Phi_p)\otimes...\otimes\Phi_{n}.
\end{equation} Since $\overline{D}^{2}=- D^{2}=-\partial_x $,
every differential operator $A\in \frak{D}_{\underline{\l},\mu}$ can be
expressed in the form (see \cite{GMO})
\begin{equation} \label{Bilinear  operators}
A=\dis\sum_{\ell=0}^{2k}\dis\sum_{|\underline{i}|=\ell}a_{\underline{i}}(x,\theta)\overline{D}^{i_{1}}\otimes
\dis\overline{D}^{i_{2}}\otimes \cdots \dis\otimes\overline{D}^{i_{n}}\end{equation} where the coefficients $a_{\underline{i}}$
are smooth functions on $S^{1|1}$ and $\ell\in\mathbb{N}$. That is, forall \\ $F_{1}=
f_{1}\alpha^{\lambda_{1}}\in\mathfrak{F}_{\l_1}, F_{2}=
f_{2}\alpha^{\lambda_{2}}\in\mathfrak{F}_{\l_2},\cdots,F_{n}=
f_{n}\alpha^{\lambda_{n}}\in\mathfrak{F}_{\l_n}$,
\begin{equation} A(F_{1}\otimes F_{2}\otimes \cdots\otimes F_{n})= \Big(\dis\sum_{\ell=0}^{2k}\dis\sum_{|\underline{i}|=\ell}a_{\underline{i}}(x,\theta)(-1)^{(\sum_{p=1}^{n-1}
|f_{p}|\sum_{s=p+1}^{n}i_{s})}\overline{D}^{i_{1}}(f_{1})
\dis\overline{D}^{i_{2}}(f_{2}) \cdots \dis\overline{D}^{i_{n}}(f_{n})\Big)\alpha^\mu.
\end{equation} Moreover, if $A\in
\mathfrak{D}_{\underline{\lambda},\mu}^k$ then $\ell=2k$. For
short, we will write the operator $A$ as:
\begin{equation}\label{operator} A=\sum_{\ell=0}^{2k}\dis\sum_{|\underline{i}|=\ell} a_{\underline{i}}\overline{D}^{\underline{i}}.\end{equation}
Where
$\overline{D}^{\underline{i}}=\overline{D}^{i_{1}}\otimes\dis\overline{D}^{i_{2}}\otimes
\cdots \otimes\dis\overline{D}^{i_{n}}$.
Thus, we consider a family of $\cK(1)$-actions on the superspace of
multilinear differential operators \\ $\frak{D}_{\underline{\l},\mu}
:=\Hom_{\rm{diff}}(\mathfrak{F}_{\l_1}\otimes\mathfrak{F}_{\l_2}\otimes
\cdots \otimes\mathfrak{F}_{\l_n},\mathfrak{F}_{\mu})$:
\begin{equation}\label{d-action}
\fL^{\underline{\l},\mu}_{X_F}(A)=\fL^{\mu}_{X_F}\circ A-(-1)^{|A||F|}
A\circ \fL^{\underline{\l}}_{X_F}.
\end{equation}

\subsection{Explicit formulas for the action of $\mathcal{K}(1)$ on
$\mathfrak{D}^{k}_{\underline{\lambda},\mu}$}
Let us calculate explicitly the action  $\mathcal{K}(1)$ on the superspace
$\mathfrak{D}^{k}_{\underline{\lambda},\mu}$. Given a differential operator
$A= \dis\sum_{\ell=0}^{2k}\dis\sum_{|\underline{i}|=\ell}a_{\underline{i}}\overline{D}^{i_{1}}
\otimes \overline{D}^{i_{2}}\otimes \cdots \otimes\overline{D}^{i_{n}}\in
\mathfrak{D}^{k}_{\underline{\lambda},\mu}$ and  ${X_F}, F\in
C^\infty(S^{1|1})$, an arbitrary contact vector field.
\begin{proposition}\label{Action}
The naturel action of $\mathcal{K}(1)$ on $\mathfrak{D}^{k}_{\underline{\lambda},\mu}$
is given by:
\begin{equation}\mathfrak{L}_{X_F}^{\underline{\lambda},\mu}(A)=
\sum_{\ell=0}^{2k}\sum_{|\underline{i}|=\ell}a^X_{\underline{i}}\overline{D}^{i_{1}} \otimes\overline{D}^{i_{2}}
\otimes \cdots \otimes\overline{D}^{i_{n}}\end{equation} where:
\begin{equation}\label{Actionformulas}\begin{array}{ll}
a_{\underline{i}}^X= \mathfrak{L}_{X_F}^{\delta-\frac{|\underline{i}|}{2}}(a_{\underline{i}}) -
{\dis\sum_{r=1}^{2k-|\underline{i}|}}(-1)^{r(|F|+|a_{\underline{i} +r\mathbf{1}_{1}}|)}\Big[\left( r+i_{1}
\atop r+2  \right)_s - \frac{1}{2}(-1)^{i_{1}}\left( r+i_{1} \atop r+1
\right)_s\\
+ \lambda_1\left( r+i_{1} \atop r  \right)_s\Big] \overline{D}^r(F')
a_{\underline{i} +r\mathbf{1}_{1}}
 - \dis\sum_{t=2}^{n}{\dis\sum_{r=1}^{2k-|\underline{i}|}}(-1)^{r(|F|+|a_{\underline{i} +r\mathbf{1}_{t}}|+i_{1}+i_{2}+ \cdots +i_{t-1})} \Big[\left( r+i_{t} \atop r+2  \right)_s - \\
\frac{1}{2}(-1)^{i_{t}}\left( r+i_{t} \atop r+1  \right)_s +
\lambda_t\left( r+i_{t} \atop r  \right)_s\Big] \overline{D}^r(F')
a_{\underline{i} +r\mathbf{1}_{t}}\end{array} \end{equation}
\end{proposition}
\begin{proof}
Let  $\phi_{1}=
\varphi_{1}\alpha^{\lambda_{1}}\in\mathfrak{F}_{\l_1}, \phi_{2}=
\varphi_{2}\alpha^{\lambda_{2}}\in\mathfrak{F}_{\l_2},\cdots,\phi_{n}=
\varphi_{n}\alpha^{\lambda_{n}}\in\mathfrak{F}_{\l_n}$. Upon using (\ref{actiondensite}), (\ref{Action on densities}) and (\ref{d-action}),  we get
\begin{equation*}
\begin{array}{ll}
   \mathfrak{L}^{\underline{\l},\mu}_{X_F}(A)\Big(\phi_{1}\otimes \phi_{2}\otimes \cdots \otimes \phi_{n}\Big) =  & \mathfrak{L}_{X_F}^\mu\Big( A(\Big(\phi_{1}\otimes \phi_{2}\otimes \cdots \otimes \phi_{n}\Big))\Big) - \\
   \cr & (-1)^{|A||F|}
A\Big(\mathfrak{L}_{X_F}^{\l_1}(\phi_1)\otimes\phi_2\otimes \cdots \otimes \phi_{n}\Big) \\
   \cr & - (-1)^{|A|(|F|+ |\phi_1|)}
A\Big(\phi_1\otimes\mathfrak{L}_{X_F}^{\l_2}(\phi_2)\otimes \cdots \otimes \phi_{n}\Big)- \cdots - \\
   \cr & (-1)^{|A|(|F|+ |\phi_1| + \cdots + |\phi_{n-1}| )}
A\Big(\phi_1\otimes \phi_2\otimes \cdots \otimes \mathfrak{L}_{X_F}^{\l_n}(\phi_{n})\Big)
\end{array}
\end{equation*}
\begin{equation*}\begin{array}{ll}
=\Big[\dis\sum_{\ell=0}^{2k}\sum_{i_{1}+i_{2}+\cdots + i_{n}=\ell}F\Big((-1)^{\sum _{j=1}^{n-1}|\varphi_{j}|(i_{j+1}+\cdots + i_{n} )} a_{\underline{i}}\overline{D}^{i_{1}}(\varphi_1)
 \overline{D}^{i_{2}}  (\varphi_2) \cdots \overline{D}^{i_{n} } (\varphi_n)\Big)'
+\\ \\
\frac{1}{2}D(F)\overline{D}\Big((-1)^{\sum _{j=1}^{n-1}|\varphi_{j}|(i_{j+1}+\cdots + i_{n} )} a_{\underline{i}}\overline{D}^{i_{1}}(\varphi_1)
 \overline{D}^{i_{2}}  (\varphi_2) \cdots \overline{D}^{i_{n} } (\varphi_n)\Big)
+ \\ \\ \mu
F'(-1)^{\sum _{j=1}^{n-1}|\varphi_{j}|(i_{j+1}+\cdots + i_{n} )} a_{\underline{i}}\overline{D}^{i_{1}}(\varphi_1)
 \overline{D}^{i_{2}}  (\varphi_2) \cdots \overline{D}^{i_{n} } (\varphi_n)
-\\ \\
(-1)^{|A||F|}(-1)^{\sum _{j=1}^{n-1}|\varphi_{j}|(i_{j+1}+\cdots + i_{n} )}
(-1)^{|F|(i_{2}+\cdots + i_{n} )} a_{\underline{i}} \\ \\ \overline{D}^{i_{1}}\Big(F\varphi_1'
+ \frac{1}{2}D(F)\overline{D}(\varphi_1)+
\l_1F'\varphi_1\Big)\overline{D}^{i_{2}}  (\varphi_2) \cdots \overline{D}^{i_{n} } (\varphi_n) - \\ \\
(-1)^{|A|(|F|+|\varphi_1|)}(-1)^{\sum _{j=1}^{n-1}|\varphi_{j}|(i_{j+1}+\cdots + i_{n} )}
(-1)^{|F|(i_{3}+\cdots + i_{n} )} \\ \\ a_{\underline{i}}\overline{D}^{i_{1}}(\varphi_1)\overline{D}^{i_{2}}\Big(F\varphi_2'
+ \frac{1}{2}D(F)\overline{D}(\varphi_2)+ \l_2F'\varphi_2 \big) \overline{D}^{i_{3}}  (\varphi_3) \cdots \overline{D}^{i_{n} } (\varphi_n) - \cdots \\ \\
 - (-1)^{|A|(|F|+|\varphi_1|+ \cdots + |\varphi_n|)}(-1)^{\sum _{j=1}^{n-1}|\varphi_{j}|(i_{j+1}+\cdots + i_{n} )}
 \\ \\ a_{\underline{i}}\overline{D}^{i_{1}}(\varphi_1)\overline{D}^{i_{2}}(\varphi_2)
 \overline{D}^{i_{3}}  (\varphi_3) \cdots \overline{D}^{i_{n} }\Big(F\varphi_n'
+ \frac{1}{2}D(F)\overline{D}(\varphi_n)+ \l_nF'\varphi_n \big)
\Big]\alpha^{\mu}.
\end{array} \end{equation*}
 Using  the super Leibniz formula (\ref{super Leibniz formula}) and by writing (\ref{Action}) in the form
\begin{equation*}\begin{array}{ll}\mathfrak{L}^{\underline{\l},\mu}_{X_F}(A)\Big(\phi_{1}\otimes \phi_{2}\otimes \cdots \otimes \phi_{n}\Big)
=\\ \Big[\dis\sum_{\ell=0}^{2k}\sum_{i_{1}+i_{2}+\cdots +
i_{n}=\ell} (-1)^{\sum _{j=1}^{n-1}|\varphi_{j}|(i_{j+1}+\cdots +
i_{n} )} a^X_{\underline{i}}\overline{D}^{i_{1}}(\varphi_1)
 \overline{D}^{i_{2}}  (\varphi_2) \cdots \overline{D}^{i_{n} } (\varphi_n) \Big]\alpha^{\mu},\end{array}\end{equation*}
By identification, we get easily the formulas (\ref{Actionformulas}).
\end{proof}

\subsection{Space of symbols of multilinear differential operators.}
Consider the graded $\cK(1)$-module
$\hbox{gr}(\mathfrak{D}^{k}_{\underline{\lambda},\mu})$ associated with the filtration
\begin{equation}\label{filtration} \mathfrak{D}_{\underline{\lambda},\mu}^0\subset\mathfrak{D}_{\underline{\lambda},\mu}^{\frac{1}{2}}
\subset\mathfrak{D}_{\underline{\lambda},\mu}^1\subset
\mathfrak{D}_{\underline{\lambda},\mu}^{\frac{3}{2}}\subset\cdots
\mathfrak{D}_{\underline{\lambda},\mu}^{k-\frac{1}{2}}\subset\mathfrak{D}_{\underline{\lambda},\mu}^{k}
\subset\cdots\end{equation}
i.e, the direct sum
\begin{equation}\hbox{gr}(\mathcal{D}_{\underline{\lambda},\mu})=
\bigoplus_{k=0}^\infty \mathfrak{D}_{\underline{\lambda},\mu}^{k} /
\mathfrak{D}_{\underline{\lambda},\mu}^{k-\frac{1}{2}}
\end{equation}
We call this $\cK(1)$-module the \textit{space of symbols of multilinear
differential operators} and denote it $\mathcal{S}_{\underline{\lambda},\mu}$. \\
The quotient module
$\mathfrak{D}^{k}_{\underline{\lambda},\mu}/\mathfrak{D}^{k-\frac{1}{2}}_{\underline{\lambda},\mu}$,
$ k\in \frac{1}{2}\mathbb{N}$, can be decomposed into   $n_{k} = \left( 2k+n-1 \atop n-1  \right)$
components that transform under coordinates change as $\delta-
\frac{k}{2}$ densities, where $\delta= \mu-|\underline{\l}|$. Therefore,
the multiplication of these components by any non-singular matrix $\varpi$
 gives rise to a $\cK(1)$-invariant isomorphism called \textit{a principal symbol
map}
\begin{equation}
\begin{CD}
\label{symbmap} \s^\varpi_{pr}:
\mathfrak{D}^{k}_{\underline{\lambda},\mu}/\mathfrak{D}^{k-\frac{1}{2}}_{\underline{\lambda},\mu}
@>\simeq>>
\mathfrak{F}_{{\delta}-\frac{k}{2}}\oplus\mathfrak{F}_{{\delta}-\frac{k}{2}}\oplus\cdots\oplus
\mathfrak{F}_{{\delta}-\frac{k}{2}}
~(n_{k} ~\rm{copies})\\
\end{CD}.
\end{equation}


 The space of symbols of order $\leq k, k\in \frac{1}{2}\mathbb{N}$, is \\
 \begin{equation}
\mathcal{S}^{k}_{ \underline{\lambda}, \mu}=\bigoplus_{\ell=0}^{2k}\mathfrak{D}^{\ell}_{\underline{\lambda},
\mu}/\mathfrak{D}^{\ell -\frac{1}{2}}_{\underline{\lambda},\mu}
\end{equation}

 The $\cK(1)$-module  $\mathcal{S}_{\underline{\lambda},\mu}$ depends only on the shift, $\d$, of the
weights and not on $\mu$, $\l_1$,  $\l_2$ $\cdots$  $\l_n$  independently. Moreover,
 for every $k \in \frac{1}{2}\mathbb{N}$, we have
\begin{equation}
\mathcal{S}^{k}_{\mu-|\underline{\lambda}|} = {\cal S}^k_{\delta}=\bigoplus_{\ell=0}^{2k}\mathfrak{D}^{\ell}_{\underline{\lambda},\mu}/
\mathfrak{D}^{\ell-\frac{1}{2}}_{\underline{\lambda},\mu}=\bigoplus_{\ell=0}^{2k}
{\mathfrak{F}}^{(\ell)}_{\delta-\frac{\ell}{2}},
\end{equation} here  the notation    ${\mathfrak{F}}^{(i)}_{\l}, i\in \mathbb{N}$ and $ \l\in
\mathbb{C}$, stands for the sum $ \bigoplus\mathfrak{F}_{\l}$ where
$\mathfrak{F}_{\l}$ is counted  $\left( i+n-1 \atop n-1  \right)$ times.\\
Thanks to the isomorphism (\ref{symbmap}), an element $P$ of ${\cal
S}^k_{\delta}$ can be written in a unique way in the form
\begin{equation}
P=\label{symboloforderk}
\a^\delta\sum_{\ell=0}^{2k}\sum_{|\underline{i}|=\ell}\bar{a}_{\underline{i}}(x,\theta)\,\a^{-\frac{|\underline{i}|}{2}}
\end{equation}
where  $\bar{a}_{\underline{i}}$ are arbitrary functions in
$C^\infty(S^{1|1})$. \\
As the orthosymplectic superalgebra $\mathfrak{osp}(1|2)$ is a subalgebra $\cK(1)$, the space of symbols $\mathcal{S}_{\delta}$ can be viewed as an $\mathfrak{osp}(1|2)$-module.

\section{ $\mathfrak{osp}(1|2)$-equivariant symbol and quantization maps.}
We restrict the $\cK(1)$-module structures to the particular subalgebra $\mathfrak{osp}(1|2)$
and look for $\mathfrak{osp}(1|2)$-isomorphisms between
$\mathfrak{D}_{\underline{\l},\mu}$ and  $\mathcal{S}_{\delta}$.
We fix a principal symbol map $\s^\varpi_{pr}$ as in (\ref{symbmap}),
where $\varpi$ is a non singular matrix.
\begin{definition}
A symbol map is a
 a linear bijection
\begin{equation}
\s^\varpi_{\underline{\l},\m}:
\mathfrak{D}_{\underline{\l},\mu}\rightarrow
\mathcal{S}_{\delta}\end{equation} such that the highest-order term
of $\s^\varpi_{\underline{\l},\m}(A)$, where $A \in
\mathfrak{D}_{\underline{\l},\mu}$,  coincides with the principal
symbol $\s^\varpi_{pr}(A)$. Hence, the inverse map,
$Q=(\s^\varpi_{\underline{\l},\m})^{-1}$, will be called a
\textit{quantization map}.
\end{definition}
 The problem of existence and uniqueness of $\mathfrak{osp}(1|2)$-equivariant symbol
(and so quantization) map can be tackled once  the symbol map
$\s^\varpi_{pr}$ is fixed.\\
The first main result of this paper is the following:
\begin{theorem} \label{MainThm}   if $\delta$ is non-resonant, i.e., $\delta=\mu-|\l|\neq
\dis\frac{1}{2}, 1,\frac{3}{2},2,\frac{5}{2},\cdots,k$ then,
$\mathfrak{D}^{k}_{\underline{\lambda},\mu}$ and $\mathcal{S}^k_{\delta}$ are
$\mathfrak{osp}(1|2)$-isomorphic through the
 family of $\mathfrak{osp}(1|2)$-equivariant maps
$\s_{\underline{\l},\mu}^\varpi$ defined by:

\begin{equation}  \s_{\underline{\l},\mu}^\varpi(A)= \alpha^\delta
\dis\sum_{p=0}^{2k}\sum_{|\underline{i}|=p}\dis\sum_{\ell=p}^{2k}
\sum_{|\underline{s}|=\ell} \varpi^{\underline{s}}_{\underline{i}} D^{\ell-p} (a_{\underline{s}})
{\alpha^{-\frac{|\underline{i}|}{2}}}
\end{equation}
where $A= \dis\sum_{p=0}^{2k}\sum_{|\underline{i}|=p} a_{\underline{i}} \overline{D}^{i_{1}} \otimes
\overline{D}^{i_{2}}\otimes...\otimes\overline{D}^{i_{n}} \in \mathfrak{D}^{k}_{\underline{\lambda},\mu} $
and $\varpi^{\underline{s}}_{\underline{i}}$ are constants  given by the induction
formula
\begin{equation}\begin{array}{ll}\label{coditioninv}
      (-1)^{\ell-p}\Big( [\frac{\ell-p}{2}]+ ( 1- (-1)^{\ell-p})(\delta -
\frac{\ell}{2})\Big) \varpi_{\underline{i}}^{\underline{s}} - \Big( [\frac{s_{1}}{2}]+ ( 1-
(-1)^{s_{1}})\l_1\Big) \varpi_{\underline{i}}^{\underline{s}-\mathbf{1}_{1}}\\ - \sum_{j=2}^{n}(-1)^{s_{1}+s_{2}+...+s_{j-1}} \Big( [\frac{s_{j}}{2}]
+ ( 1- (-1)^{s_{j}})\l_j\Big)\varpi_{\underline{i}}^{\underline{s}-\mathbf{1}_{j}}=0.
     \end{array} \end{equation}
If $\varpi$ is the idendity map, we obtain the "normalized" symbol map $\s_{\underline{\l},\mu}^{Id}$  given by the rule
\begin{equation}\s_{\underline{\l},\mu}^{Id}(A)= \alpha^\delta
\dis\sum_{p=0}^{2k}\sum_{|\underline{i}|=p}\dis\sum_{\ell=p}^{2k}
\sum_{|\underline{s}|=\ell\atop s_{1}\geq i_{1},s_{2}\geq i_{2},...,s_{n}\geq i_{n}}\gamma^{\underline{s}}_{\underline{i}}D^{\ell-p}
(a_{\underline{s}}) {\alpha^{-\frac{|\underline{i}|}{2}}}\end{equation} such that

 \begin{equation} \label{SymMapMainFor}  \gamma^{\underline{s}}_{\underline{i}}=
(-1)^{[\frac{\ell-p+1}{2}]}\prod_{t=2}^{n}\dis\frac{
(-1)^{\psi(t)}\left( \varphi(t) \atop s_{t}-i_{t} \right)_{s}
\Xi_{s_{t},i_{t}}(\lambda_{t})} {\left( [\frac{\varphi(t)}{2}] \atop
[\frac{s_{t}-i_{t}}{2}] \right) \left( [\frac{\varphi(t)+1}{2}]
\atop [\frac{s_{t}-i_{t}+1}{2}] \right) \left(
2\d-p-1\atop  [\frac{\ell-p+1}{2}] \right)} \Xi_{s_{1},i_{1}}(\lambda_{1})\end{equation}\\
where the functions $\varphi, \psi$ and $\Xi$ are defined by\\
$\varphi(t)=\dis\sum_{j=1}^{t}s_{t}-i_{t}$,\ $\psi(t)=\dis\sum_{j=1}^{t-1}s_{j}(s_{j+1}-i_{j+1})$,\  $\Xi_{s_{t},i_{t}}(\l_{t})= \left([\frac{s_{t}}{2}]\atop [\frac{i_{t}}{2}] \right)\left( 2\l_{t}+[\frac{s_{t}-1}{2}] \atop [\frac{2(s_{t}-i_{t})+1+(-1)^{i_{t}}}{4}]
\right)$\\
and the notation $\left( \nu\atop q \right)$ stands for the binomial coefficient given  by
$\left( \nu\atop q \right)=\frac{\nu(\nu-1)\cdots(\nu-q+1)}{q!}$.\\
Moreover, once the principal symbol is fixed, the symbol map
$\s_{\underline{\l},\mu}^\varpi$ is unique.
\end{theorem}
\begin{proof}
We begin the proof by proving the $\mathfrak{osp}(1|2)$-equivariance
of the map $\s_{\underline{\l},\mu}^{Id}$. Indeed, Let
$X=X_F\in\mathcal{K}(1)$. We have\\
$\s_{\underline{\l},\mu}^{Id}\Big((\mathfrak{L}_{X}^{\underline{\l},\mu}(
A)\Big)= \alpha^\delta
\dis\sum_{p=0}^{2k}\sum_{|\underline{i}|=p}\overline{a}^{X}_{\underline{i}}{\alpha^{-\frac{|\underline{i}|}{2}}}$.
Then, we readily see that
\begin{equation*}\overline{a}^X_{\underline{i}}=\sum_{\ell=p}^{2k}\sum_{ |\underline{s}|=\ell}\gamma_{\underline{i}}^{\underline{s}}D^{(\ell-p)}({a}^X_{\underline{i}}),\ p=|\underline{i}|.\end{equation*}
Thanks to the proposition \ref{Action}, for all $0\leq p=|\underline{i}|\leq k$,
we get\\
\begin{equation*}\begin{array}{ll}&\overline{a}^X_{\underline{i}}=\dis\sum_{\ell=p}^{2k}\dis\sum_{ |\underline{s}|=\ell}(-1)^{|a^X_{\underline{s}}|(\ell-p) + \frac{(\ell-p)(\ell-p+1)}{2}}\gamma_{\underline{i}}^{\underline{s}}\overline{D}^{(\ell-p)}({a}^X_{\underline{s}}) \\
&=  \dis\sum_{\ell=p}^{2k}\dis\sum_{
|\underline{s}|=\ell}(-1)^{|a^X_{\underline{s}}|(\ell-p) +
\frac{(\ell-p)(\ell-p+1)}{2}}\gamma_{\underline{i}}^{\underline{s}}\overline{D}^{(\ell-p)}
\Big[\mathfrak{L}_{X_F}^{\delta-\frac{|\underline{s}|}{2}}(a_{\underline{s}})   \\
&-{\dis\sum_{m=1}^{2k-|\underline{s}|}}(-1)^{m(|F|+|a_{\underline{s}+\mathbf{1}_{1}}|)}\Big(\left(
m+s_{1} \atop m+2  \right)_s - \frac{1}{2}(-1)^{s}\left( m+s_{1} \atop m+1
\right)_s + \lambda_1\left( m+s_{1} \atop m  \right)_s\Big) \dis\overline{D}^m(F') a_{\underline{s}+m\mathbf{1}_{1}}  \\
 &-{\dis\sum_{m=1}^{2k-|\underline{s}|}}\sum_{j=2}^{n}(-1)^{m(|F|+|a_{\underline{s}+m\mathbf{1}_{j}}|+
 \sum_{t=1}^{j-1}s_{t})} \\
  & \Big(\left( m+s_{j} \atop m+2  \right)_s -
\frac{1}{2}(-1)^{s_{j}}\left( m+s_{j} \atop m+1  \right)_s +
\lambda_j\left( m+s_{j} \atop m \right)_s\Big)
\dis\overline{D}^m(F') a_{\underline{s}+m\mathbf{1}_{j}} \Big].
\end{array}\end{equation*}
Thus
\begin{equation*}\begin{array}{ll}
\overline{a}^X_{\underline{i}}-
\mathfrak{L}_X^{(\d-\frac{p}{2})}(\overline{a}_{\underline{i}})=
{\dis\sum_{\ell=p}^{2k}\sum_{
|\underline{s}|=\ell}}(-1)^{|a_{\underline{s}}|+|F|+(\ell-p)}\gamma_{\underline{i}}^{\underline{s}}
\Big[(\delta -\frac{l}{2})\left( \ell-p\atop 1\right)_s +
\frac{1}{2}(-1)^{\ell-p}\left( \ell-p \atop 2\right)_s \\+ \left(
\ell-p\atop 3\right)_s \Big] \overline{D}(F') D^{\ell-p-1}( a_{\underline{i}})
- {\dis\sum_{\ell=p}^{2k}\sum_{ |\underline{s}|=\ell}}(-1)^{|a_{\underline{s}}|+|F| }\Big[
\gamma_{\underline{i}}^{\underline{s}-\mathbf{1}_{1}} \Big( \left( s_{1}\atop 3  \right)_s +
\frac{1}{2}(-1)^{s_{1}}\left( s_{1} \atop 2 \right)_s +    \l_1 \left(
s_{1} \atop 3  \right)_s \Big) \\+ \dis\sum_{j=2}^{n}(-1)^{\sum_{t=1}^{j-1}s_{t}} \gamma_{\underline{i}}^{\underline{s}-\mathbf{1}_{j}} \Big(
\left( s_{j}\atop 3  \right)_s + \frac{1}{2}(-1)^{s_{j}}\left( s_{j} \atop 2
\right)_s +   \l_{j} \left( s_{j} \atop 3  \right)_s  \Big)
\Big]\overline{D}(F') D^{\ell-p-1}( a_{\underline{i}}) + \\ ~(\hbox{higher terms
in}~\overline{D}^{n}(F'), n\geq2).
\end{array}\end{equation*}
Now, through a simple calculation, one can check out that the
scalars $\gamma_{\underline{i}}^{\underline{s}}$ satisfis the relationship
\begin{equation*}
      (-1)^{\ell-p} \Upsilon(\delta - \frac{\ell}{2},\ell-p) \gamma_{\underline{i}}^{\underline{s}} -   \Upsilon(\l_1,s_{1}) \gamma_{\underline{i}}^{\underline{s}-\mathbf{1}_{1}} - \sum_{j=2}^{n}(-1)^{s_1+s_2+...+s_{j-1}} \Upsilon(\l_j,s_{j})\gamma_{\underline{i}}^{\underline{s}-\mathbf{1}_{j}}=0
      \end{equation*} where, for $\l\in \mathbb{C}$ and $ m\in \mathbb{N}$, we put
$$\Upsilon(\l,m) = \frac{1}{2}\big( [\frac{m}{2}]+ ( 1-
(-1)^m)\l\big).$$
Since, the term in $\overline{D}(F')$ vanishes, we can clearly see that the map $\s_{\underline{\l},\mu}^{Id}$  $\mathfrak{osp}(1|2)$-equivariant .\\
Now, we can easily adapt the
proof of locality given in \cite{GMO} for the unary case to our
case and then use the locality property of an
$\mathfrak{osp}(1|2)$-equivariant symbol map. Therefore, in addition, from the expression of the "normalized" symbol
map $\s_{\underline{\l},\mu}^{Id}$ we can suppose that a general symbol map
$\s_{\underline{\l},\mu}^\varpi$ can be written as

\begin{equation} \label{generalsym} A= \dis\sum_{p=0}^{2k}\sum_{|\underline{i}|=p} a_{|\underline{i}|} \overline{D}^{i_1} \otimes
\overline{D}^{i_2}\otimes...\otimes\overline{D}^{i_n} \longmapsto \alpha^\delta
\dis\sum_{p=0}^{2k}\sum_{|\underline{i}|=p}\dis\sum_{\ell=p}^{2k}
\sum_{|\underline{s}|=\ell} \varpi^{\underline{s}}_{\underline{i}}(x,\theta) D^{\ell-p} (a_{\underline{s}})
{\alpha^{-\frac{|\underline{i}|}{2}}}.
\end{equation}
Obviously, to get the condition of
$\mathfrak{osp}(1|2)$-equivariance, it is sufficent to impose invariance
with respect to the vector fields $D= 2X_\theta$ and $xD=
2X_{x\theta}$ to meet the whole condition $\mathfrak{osp}(1|2)$-equivariance. Thus we have:
\begin{description}
  \item[a)]  A symbol map
(\ref{generalsym}) commutes with the action of $D $ if and only if
the coefficients $\varpi^{\underline{s}}_{\underline{i}}$  are
constants (i.e., do not depend on $x, \theta$),
  \item[b)]  A symbol map
(\ref{generalsym}) commutes with the action of $xD $ if and only if the coefficients
$\varpi^{\underline{s}}_{\underline{i}}$ satisfy
the induction formula (\ref{coditioninv}).
\end{description}
If $\delta= \mu-|\underline{\l}|$ is non-resonant, i.e.,
$\delta=\mu-|\underline{\l}|\neq \dis\frac{1}{2},
1,\frac{3}{2},2,\frac{5}{2},\cdots,k$ , then, it is easy to see that
the solution of the equation (\ref{coditioninv}) and once the
principal symbol $\s^\varpi$ where
$\varpi=(w_{\underline{i}}^{\underline{i}})_{ |\underline{i}|=2k}$
is fixed, the symbol map $\s_{\underline{\l},\mu}^\varpi$ is unique.
\end{proof}
\begin{remark} We can write the symbol map $\s_{\underline{\l},\mu}^{Id}$ as in \cite{GMO}, Theorem $6.1$. Indeed
Let $A=a_{\underline{i}}\overline{D}^{i_{1}} \otimes \overline{D}^{i_{2}}\otimes...\otimes\overline{D}^{i_{n}}\in
\mathfrak{D}^{k}_{\underline{\lambda},\mu}$ and $|\underline{i}|=2k$, then

\begin{equation} \s^{Id}_{\underline{\l},\mu}(A)= \alpha^\delta
\dis\sum_{\ell=0}^{2k} \sum_{|\underline{s}|=\ell} \chi_{\underline{i}}^{\underline{s}} D^\ell
(a_{\underline{i}}) {\alpha^{\frac{|\underline{s}|-|\underline{i}|}{2}}}
\end{equation}
where

\begin{equation}  \chi_{\underline{s}}^{\underline{i}}=\left\{\begin{array}{l}
(-1)^{[\frac{\ell+1}{2}]}\dis\prod_{t=2}^{n}\dis\frac{(-1)^{\Delta(t)}
\left( \Gamma(t) \atop s_{t} \right)_{s}
\Xi_{i_{t},i_{t}-s_{t}}(\l_{t}) } {\left( [\frac{\Gamma(t)}{2}]
\atop [\frac{s_{t}}{2}] \right) \left( [\frac{\Gamma(t)+1}{2}] \atop
[\frac{s_{t}+1}{2}] \right) \left( 2\d+\varphi(n)-1\atop
[\frac{\ell+1}{2}]
\right)}\Xi_{i_{1},i_{1}-s_{1}}(\l_{1})\ \\ \hbox{if}\ i_{t}\geq s_{t},  t \in \{1,2,...,n\}, \\
 0 \ \hbox{otherwhise}\end{array}\right.,\end{equation}
 here
 $\Gamma(t)=\dis\sum_{j=1}^{t}s_{j}$ and $\Delta(t)=\dis\sum_{j=1}^{t-1}s_{j}s_{j+1}$.
\end{remark}

Now, by a direct computation, one can easily check  the following  explicit formula for the quantization map
$Q_{\underline{\l},\mu}^{Id}$:
\begin{proposition}
\label{MainThmbis} The quantization map $Q_{\underline{\l},\mu}^{Id}$, i.e.,
the inverse of the symbol map $\s^{Id}_{\underline{\l},\mu}$ given in
theorem \ref{MainThm} associates to a polynomial $P= \alpha^\delta
\dis\sum_{\ell=0}^{2k}\sum_{|\underline{i}|=\ell}\bar{b}_{\underline{i}}{\alpha^{-\frac{|\underline{i}|}{2}}}\in
\mathcal{S}_{\delta}^k $ the  differential operator
$Q_{\underline{\l},\mu}^{Id}(P)=
\dis\sum_{\ell=0}^{2k}\sum_{|\underline{i}|=p}\tilde{b}_{\underline{i}}\overline{D}^{i_{1}}\otimes\overline{D}^{i_{2}}\otimes...\otimes\overline{D}^{i_{n}}
\in \mathfrak{D}^{k}_{\underline{\lambda},\mu}$ such that\\
$\tilde{b}_{\underline{i}}=\dis\sum_{\ell=p}^{2k}\sum_{\underline{s}=\ell}\beta_{\underline{i}}^{\underline{s}}D^{\ell-p}(\bar{b}_{\underline{s}}),$
where

\begin{equation}\left\{\begin{array}{ll}
\beta_{\underline{i}}^{\underline{s}}=
(-1)^{[\frac{\ell-p-1}{2}]}\dis\prod_{t=2}^{n}\dis\frac{
(-1)^{\psi(t)}\left( \varphi(t) \atop s_{t}-i_{t} \right)_{s}
\Xi_{s_{t},i_{t}}(\lambda_{t})} {\left( [\frac{\varphi(t)}{2}] \atop
[\frac{s_{t}-i_{t}}{2}] \right) \left( [\frac{\varphi(t)+1}{2}]
\atop [\frac{s_{t}-i_{t}+1}{2}] \right) \left(
2\d-l\atop  [\frac{\ell-p+1}{2}] \right)}\Xi_{s_{1},i_{1}}(\lambda_{1})\\
  ~ \hbox{if}~ \ell=|\underline{s}|> p=|\underline{i}|\\
\beta_{\underline{i}}^{\underline{s}}=
\gamma_{\underline{i}}^{\underline{s}} ~~~~~~~ \hbox{if}~~~
|\underline{s}|=|\underline{i}|
\end{array}\right.\end{equation}.
\end{proposition}

\end{document}